\documentclass{amsart}

\title{AP Theory III: Cone-like Graded SUSY, Dynamic Dark Energy and the YM Millenium Problem}
\author{H. E. Winkelnkemper}

\address{\begin{flushleft}\quad Department of Mathematics \\
\quad University of Maryland \\
\quad College Park, Maryland 20742
\end{flushleft}}

\email{hew@math.umd.edu}

\begin{document}
\maketitle

\begin{abstract} Artin Presentation Theory, (AP Theory), is a new, direct infusion, via pure braid theory, of discrete group theory, (i.e.,  symmetry in its purest form), into the theory of {\it smooth} 4-manifolds, (i.e., $(3+1)$-Quantum Gravity in its purest topological form), thus exhibiting the most basic, rigorous, universal, model-free intrinsic {\it gauge-gravity duality}, in a non-infinitesimal, cone-like graded, as holographic as possible,  model-independent, non-perturbative, background-independent, parameter-free manner. {\it In AP Theory even smooth topology change becomes gauge-theoretic, thus setting the stage for a rigorous smooth topological $(3+1)$-QFT of Dynamic Dark Energy}. In this theory, the rigid $\infty$ of the dimension of classical Hilbert space is substituted by the dynamic $\infty$ of the $\infty$ generation at each stage of a cone-like graded group of topology-changing transitions/interactions. As a corollary, the Cosmological Constant problem and the YM Millenium Mass Gap problem, two of the most perplexing main problems of modern theoretical physics, become rigorously, intimately mathematically related, by having the same qualitative {\it dynamical} roots. Ultimately our main point is meta-mathematical, as far as modern physics is concerned: due to the discrete group-theoretic conceptual simplicity of the theory, with its group-theoretic 'Planckian membrane/discreteness' starting point, {\it the fact that it is not just a mere mathematical model}, and all its  properties above, any other {\it mathematically rigorous} approach to the above problems has to built on AP Theory and be topologically  absorbed and enveloped by it.

\end{abstract}

\section{Introduction}

{\it And now we might add something concerning a certain most subtle Spirit, which pervades and lies hid in all gross bodies.}
I. Newton, (as quoted by Schwinger, \cite{S}).

{\it Ou bien il n'y aurait rien au monde qui ne f\^ut d'origine \'electromagn\'etique.
 
Ou bien cette partie qui serait pour ainsi dire commune \'a tous le ph\'enomenes physiques ne serait qu'une apparence, quelque chose qui tendrait \'a nos methodes de mesure.}
H. Poincar\'e, \cite{P2}, p.498.

{\it Die Frage, was das elektromagnetische Feld ist, kann so wenig beantwortet werden wie die Frage, was Masse ist.}
H. Weyl, \cite{We}, p.307.

{\it I think however, one should keep in mind that there is some unknown field involved here.}
W. Pauli, \cite{P}, p.xviii; letter to V. Weisskopf after the discovery of parity violation.

This paper exhibits a sporadic, characteristic of dimension $4$, mathematically rigorous, {\it smooth} $(3+1)$-dimensional topological theory, in which many very important analytic and quantitative problems of modern physics, some not even yet rigorously well-posed, have actual rigorous, {\it qualitative} smooth topological analogues with actual {\it dynamic} smooth topological solutions.

The Theory of Artin Presentations is a new, rigorous, direct infusion of discrete group theory, via pure framed/colored braid theory and their defining purely group theoretic equation, the Artin Equation, into the very physical relevant theory of {\it smooth}, compact $(3+1)$-manifolds, which goes beyond, (as explained on p.9 of \cite{W2}), merely substituting a 'surgery diagram' by a framed braid.

For a quick first introduction to the basic mathematics of AP Theory, all of whose concepts were already very well-known to Poincar\'e, see \cite{C1} and section 2 of \cite{W2}, (which we assume the reader has at hand), and references therein.

This new infusion of discrete group theory is obtained with group {\it Presentation Theory}, i.e., Geometric Group Theory, introduced by Poincar\'e into Topology, which is more basic and canonical than group {\it Representation Theory} of Lie groups and their subgroups, which is more familiar and well-known among physicists.

In fact, AP Theory, although very conceptually simple and absolutely mathematically rigorous, is somewhat unconventional and subtle, as far as Physics is concerned: its basic equation, the Artin Equation, which characterizes pure framed braids, is at each stage, $n$, a discrete, purely group theoretic equation in the free group, $F_n$, thus avoiding all infinitesimal analytic differential equations and their analytic difficulties and nomographic restrictions. A priori, so-called 'field theory divergences' do not appear in AP Theory; compare to \cite{De}.

 In AP Theory all its smooth differential 3-manifolds and 4-manifolds, $M^3(r)$, $W^4(r)$, are obtained by means of a  non-infinitesimal, non-local, holographic, topological construction, (\cite{W1}, \cite{W2}, \cite{C}, \cite{CW}),
starting from a diffeomorphism $h(r):\Omega_n\to \Omega_n$, determined by the Artin presentation $r$,  which restricts to the identity on the boundary of $\Omega_n$, the compact 2-disk in the plane with $n$ holes.

The Artin Equation is intimately related to Donaldson's Theorem, \cite{W2}, p.9, despite the absence of moduli in AP Theory.

In the cone-like graded AP Theory, the purely group-theoretic Artin presentation $r$, on $n$ generators, {\it a completely discrete concept}, uniquely determines, up to isotopies keeping boundaries fixed, the diffeomorphism $h(r):\Omega_n\to \Omega_n$, {\it a smooth $2D$, continuous, i.e., 'membranic' concept}.

{\bf This is the rigorous group-theoretic, {\it qualitative} topological Planckian analogue in AP Theory to Planck's {\it quantitative} analytical $E=h\nu$ {\it postulate,} i.e., Planck's Law.}

It was discovered in 1975, by Gonz\'alez-Acu\~na using work of Artin; see \cite{W1}.

Observe that Planck's Law was a postulate, whereas Gonz\'alez-Acu\~na's  topological analogue is a theorem, see \cite{W1}.
 
{\it In AP Theory the particle-wave duality becomes particle-membrane duality.}. In fact, {\it string-membrane duality}, if we consider the discrete Artin presentation $r$ intuitively as a 'cosmic string', instead of a discrete particle, see section 3 ahead.

AP Theory is the ultimate, most basic, cone-like graded Membrane Theory, emerging out of the discrete 'vacuum' of discrete group theory.

In this sense, AP Theory is a rigorous Membranic Higgs Mechanism, \cite{W}, p.96.

In AP Theory, 'mass' is this membranic mass, enhanced by the dynamic energy of the Torelli actions discussed below; compare to \cite{W}, p.132.

Topologically, 'electro-magnetism' seems to become 'membranism', a fact which appears to have been intuited in the above quotes.

It also seems to be relevant to Witten's {\it 'big question'} in \cite{Wi3}, related to the LHC: {\it 'This is the question of why "electromagnetism" is so different from the "weak interactions".'} Compare also to \cite{NN}.

All this is genuinely in the spirit of Poincar\'e 's conception of quantum physics, see remark 1 in the last section.

{\it This is the basic, fundamental, topological 'Planckian' continuum/discreteness  Duality of AP Theory}, which  also induces immediately a gauge/ gravity duality, by means of the holographic construction of smooth $(3+1)$-manifolds $W^4(r)$, from the $h(r):\Omega_n\to \Omega_n$, see \cite{W1}.

(Notice that this construction goes in exactly the opposite direction of the still heuristic $4D$ dimensional reduction to $2D$ via $\sigma$ models as in, e.g., \cite{BJSV}, as well as that at the beginning of the Geometric Langlands program, \cite{Wi5}.)

These properties allow AP Theory to have a graded, (by the positive integers), {\it infinitely generated} at each stage, group of homology preserving, but {\it smooth topology-changing} group of transition and interactions on smooth $(3+1)$-manifolds, the Torelli transitions/interactions (see \cite{W2}, \cite{C}). 

In fact, these transitions/interactions should be regarded as providing the 'integration' of the discrete Artin Equation.

This immediately supports C. N. Yang's fundamental gauge principle that "Symmetry dictates Interaction",(see \cite{D}, p. 218, also \cite{Fe}).

These Torelli transitions/interactions act naturally, simultaneously, holographically, {\it in unison}, on the Artin presentation $r$, the associated membrane $h(r):\Omega_n\to \Omega_n$ and the associated smooth $(3+1)$-manifolds $W^4(r)$; thus the  Torelli also act on the constituitivity of the space-time universes, the $W^4(r)$. Compare to Schwinger's Newton quote above.

AP Theory is, a priori, a {\it universal, intrinsic, model-independent} gauge theory, not restricted to the 'interior' metric connection gauge theory of one single given 4-manifold, just as Cobordism Theory is a universal, intrinsic model-free homology theory, not restricted to the interior homology of one single, given manifold. 

AP Theory is not related to any particular Lie group, except $SU(2)$, and that is via Donaldson's theorem, see \cite{W2}, p.9., not, a priori, with any direct $SU(2)$-connection related Differential Geometry.

Discrete group theory acts, in general, exteriorly, universally on the smooth $(3+1)$-manifolds, but, as in  Cobordism theory, it can also act on the same underlying topological $(3+1)$-manifold, (see \cite{C}), thus acting 'internally' on its constituitivity, justifyfing again the name topological gauge theory.

The discrete $r$ are the AP theoretic  'cosmic strings', 'magnetons', 'gauge particles', \cite{MO}; the $h(r):\Omega_n\to \Omega_n$, the 'membratons', 'strangelets', 'wavicles', (\cite{W}, p.132) and the $W^4(r)$, the 'space-time universes'.

{\bf The infinity of the {\it static} dimension of classical Hilbert space has been substituted by the cone-like, graded {\it dynamic} infinity, at each stage, of the generation of the cone-like graded group of homology-preserving, topology-changing Torelli transitions/interactions.}

This increases the power of AP Theory as a mere $2D$ membrane theory and will reveal the rigidity of the $\infty$ dimension of classical Hilbert space as a serious obstruction to solving important problems of Modern Physics rigorously, by considering their smooth topological 'completions' and dynamics therein, in AP Theory. 

In AP Theory this Rigidity has been broken, {\it quantized itself}, so to speak, thus also affecting classical SUSY, compare to \cite{Wi2} and section 3 ahead, as well as the classical formulations of the Cosmological Constant, \cite{Ca}, \cite{Ca2}, and Mass Gap problems, \cite{JW}, \cite{V}, see sections 4 and 5 ahead.

 Non-commutativity  appears naturally in AP Theory: the above graded group of topology-changing transitions is at each stage, $n$, isomorphic to the infinitely generated commutator subgroup of the pure braid group on $n$ strands, (see \cite{W1}, p.250), with its recently discovered 'binary free-ness' property, \cite{LM}.

 Indeterminism also appears in AP Theory, but only at the most basic $2D$ 'Planckian' level: the $h(r):\Omega_n\to \Omega_n$ are only determined by the Artin presentation $r$ up to smooth isotopies, which restrict to the identity on the boundary of $\Omega_n$.

Notice here, that the vacuum fluctuations/excitations, $h(r):\Omega_n\to \Omega_n$, define 'indeterminism' rigorously, not viceversa and non-rigorously, by appealing to 'virtual' particles via the Heisenberg Uncertainty principle. Compare to \cite{H}.

This is all that is left of Bohr's {\it 'statistic correspondence'} implied by the classical quantum mechanical discontinuum/continuum correspondence, compare to \cite{P}, p.1086, and \cite{H}.

As in the work of Poincar\'e cited in section 6, remark 1, their is no need, a priori, to consider thermodynamic arguments in AP Theory.

In AP Theory this continuum/discreteness correspondence, becomes a concise cone-like, rigorous,  group-theoretic discrete/membranic correspondence, which, furthermore, is intimately related to the important concept of Quantum Chromodynamics, (QCD): 'confinement', p.130, \cite{W}; \cite{Wi4}, p.2  and section 5, ahead.

Confinement is a {\it rigorous} topological concept in AP Theory; compare to \cite{Wi4}, p.2. 

{\it The indeterminism in AP Theory immediately causes confinement.} See section 5.

 These 3 fundamental 'quantum' properties, (the $\infty$ of infinite generation, non-commutativity, indeterminism), together with the above topological Planckian membranic continuum/discreteness duality {\it as a rigorous mathematical starting point}, the fact that {\it AP Theory is not a mere mathematical model} and its discrete group-theoretic conceptual simplicity, certainly makes AP Theory into a genuine, autonomous, cone-like, graded $(3+1)$-Quantum  Gauge Theory, related moreover to the Lie group $SU(2)$, via Donaldson's Theorem. 

In the following sections, we show how AP Theory has genuine {\it qualitative} dynamic topological analogues to the three very important {\it quantitative} physical problems above and how, furthermore,  by starting with the $r$ and their corresponding $h(r):\Omega_n\to \Omega_n$ as its rigorous smooth topological Planckian starting point, its {\it 'zero point vacuum fluctuations'}, AP Theory also gives their rigorous, {\it dynamic} smooth {\it topological} solutions.

It will become clear that these three important problems can be solved in AP Theory, only because its basic equation, the Artin Equation, is not a nomographically challenged analytic, nor infinite dimensional Hilbert space using equation, such as the Schroedinger equation, the Wheeler-DeWitt equation, analytic Lagrangians, Hamiltonians, Aharonov-Bohm potentials,etc.

{\bf The Artin Equation allows us to intimately relate the most global $(3+1)$-topology-changing transition/interactions to a most local, smooth $2D$ topological Planckian starting point of discrete/continuum duality.}

In AP Theory, the most {\it local} discrete particles, the cosmic strings $r$, are related to the most {\it global} smooth universes $W^4(r)$, both tied  and bound together by the powerful symmetries, the Torelli transitions, acting {\it in unison} on both.

This is mathematically possible only because, in AP Theory, we start from a rigorous, natural topological analogue to Planck's Law and end, not with the {\it static} rigidity of infinite dimensional Hilbert space of classic Heisenberg-Schroedinger Quantum Mechanics, but with the {\it dynamic} infinite generation at each stage, of the above graded group of topology-changing transitions and interactions on the smooth $(3+1)$-spacetime universes $W^4(r)$.

This is especially relevant for topologically, qualitatively {\it dynamically} relating and solving the Cosmological Constant and YM Millenium Mass Gap problems of sections 4 and 5 below, respectively.

Since the smooth $(3+1)$-spacetime universes $W^4(r)$ represent the last topological vestige of Einsteinian Gravity, this makes AP Theory also a genuine topological Quantum Gravity Theory, a smooth topological $(3+1)$-TOE,('theory of everything'), rigorously avoiding so-called field theory divergences, black holes, etc. Compare to \cite{De}.

We conjecture that AP Theory is the unique, rigorous, 'tightest', most conceptually simple, universal metamathematical 'net' to have both these fundamental features from Quantum Mechanics and General Relativity.

It is the extraordinary conceptually simplicity of the cone-like graded symmetry of AP Theory, with its powerful Torelli transitions/interactions, which allows it to topologically absorb and mathematically relate many general, important, unsolved questions of Modern Physics to each other.

\section{AP Theory as the $(3+1)$-QFT with the Most General and Strongest Interactions}

AP Theory is, purely mathematically speaking, first of all, the sharpest topological $(3+1)$-QFT, extending this format of Atiyah, et al., to the maximum, by basing it rigorously on a 'topological Planckian Continuous to Discrete' starting point, from which this whole sporadic $(3+1)$-QFT, characteristic of dimension $4D$, follows.

Thus the 'cut and paste' categorical cobordisms, partition functions and Morse theories of classical $(3+1)$-TQFTs, are substituted by the much more subtle, graded, infinitely generated at each stage, group of topology-changing, but homology preserving group of Torelli transitions/interactions. See, e.g., the example on p.8 of \cite{W2}, as well as \cite{C}, where it is shown that such transitions can change the smooth structure of a smooth $(3+1)$-manifold, {\it but keep the underlying topological structure intact.}

{\it In AP Theory, even smooth topology-change, e.g., Morse Theory, becomes purely group-theoretic, i.e. gauge-theoretic in a model-free manner.}

The last vestiges of Hilbert space in AP Theory, as a classic $(3+1)$-TQFT, are the symmetric, integer matrices $A(r)$, the 'exponent sum' matrices associated to each Artin presentation $r$, see \cite{W3}, which are intimately related to Donaldson's Theorem, see \cite{W2}, p.9. 

AP Theory is a sporadic, holographic theory, characteristic of smooth dimension $(3+1)$, which, nevertheless, contains {\it all} closed, orientable $3$-manifolds and enough {\it smooth} $4$-manifolds, (including the Kummer surface and all simply-connected, elliptic $E(n)$ complex surfaces, \cite{CW}), that  non-trivial, purely group-theoretic Donaldson theories and Seiberg-Witten theories are still present in AP Theory (see p.9, \cite{W2}), unlike as with the $(3+1)$-TQFTs studied in, e.g., \cite{CFW}, p.2.

AP Theory still contains an {\it unavoidable},  intrinsic, irreducible concept of {\it indeterminism}: the purely discrete Artin presentation $r$, only determines the analytic, smooth $2D$ 'membranic' vacuum fluctuation/excitation $h(r):\Omega_n\to \Omega_n$ up to smooth isotopies of $\Omega_n$, which restrict to the identity on the boundary of $\Omega_n$. As mentioned above, this leads to the analogue in AP Theory of QCD 'confinement', see section 5 ahead.

Except for this smooth topological uncertainty, AP Theory is totally deterministic, compare to \cite{H}. Indeterminacy, uncertaintity in AP Theory has, so to speak, been 'quantized' itself: the vacuum fluctuation/excitations $h(r):\Omega_n\to \Omega_n$, give indeterminism, not viceversa, as with classical Heisenberg Uncertainty defining 'self-annhilating virtual particles', etc., in a non-rigorous heuristic manner.

Thus, it is reasonable to call this $(3+1)$-QFT, with the most radical topology-changing interactions {\it The $(3+1)$-QFT for Dynamic Dark Energy}, see section 4 ahead.  

What other analogues to Dynamical Dark Energy can exist in rigorous mathematical physics, than the most radical smooth $(3+1)$-topology-change possible, which moreover is directly related to a genuine topological Planckian starting point? 

If membranicity is mass,  AP Theory is {\it 'an unseen quantum field that suffuses the entire cosmos...and imparts mass to all particles'}. Compare to \cite{A}.

It other words AP Theory is the 'membranic Higgs field'. Compare to \cite{W}, p.96.

AP Theory is indeed the most conceptually simple and most universal rigorous $(3+1)$-QFT; compare to \cite{ACK}, which furthermore, is graded and cone-like, as holographic as possible, model-independent, non-perturbative, back-ground independent and parameter-free, all very desired important properties of Modern Physics.

Thus in AP Theory, at least at the smooth topological level, indeed  $(3+1)$-QFT {\it "has been developed as a mathematical subject"}, (\cite{JW}, p.4), at least smooth topologically, with dynamic applications, {\it in unison} to the Cosmological Constant problem and the YM Millenium Mass Gap problem, see sections 4 and 5 ahead.

\section{AP Theory as The Cone-like, Graded Embodiment of SUSY}

Due to its extremely basic discrete group-theoretic nature, which is graded by the positive integers, i.e. AP Theory is 'cone-like', which allows the existence of the $\infty$-generated at each stage, graded, radical topology-changing interactions in the previous section, (dynamic dark energy), it is natural to consider the sheer existence of AP Theory already as an incarnation, embodiment, so to speak, of Supersymmetry (SUSY). The sheer mathematical existence of the rigorous group-theoretic AP Theory meta-mathematically takes the place of {\it "SUSY as a lost symmetry that existed in the early universe."} See \cite{A}.

The sheer mathematcal existence of AP Theory as a whole {\it is} its rigorous analog of the "Big Bang", when all that existed was 'symmetry' balled into a 'superforce', whose analog in AP Theory are the $\infty$-generated at each stage, Torelli transitions/interactions; the Artin presentations $r$ are the associated 'cosmic strings'.

We quote Wilczek on Quantum Chromodynamics, \cite{W}, p.196: {\it "QCD is in a profound and literal sense constructed as the embodiment of symmetry".}

Due to its basic, discrete group-theoretic nature, not only is AP Theory even more universally symmetric, but, unlike QCD at this stage, it is also a {\it entirely mathematically rigorous} theory.

{\it AP Theory is the mathematically rigorous embodiment of cone-like graded SUSY.}

In AP Theory, the holographic construction of $4D$ smooth structures on the $W^4(r)$ from the $2D$ membrane diffeomorphisms, the $h(r):\Omega_n\to \Omega_n$, allowed for the existence of the radical, powerful symmetries, the Torelli, to act in unison, simultaneously, on the discrete $r$, the $h(r):\Omega_n\to \Omega_n$ and the $W^4(r)$. This is a rigorous version of Wilczek's {\it "Eliminating mass enables us to bring more symmetry into the mathematical description of Nature"}, \cite{Wil2}, p. 35. 

There is no need to evoke infinite dimensional Hilbert space for establishing a SUSY: the infinite generation, at each stage, of the action of the topology-changing Torelli transitions/interactions, represent a SUSY in a more direct, stronger {\it dynamical}, i.e., physical, manner. 

SUSY is 'broken' mainly because, in AP Theory, the rigidity of the $\infty$ dimension of classical Hilbert space has been 'broken'  in a very canonical graded, cone-like, rigorous way. Compare to \cite{Wi2}.

Other approaches to SUSY also exist in AP Theory, since one can stabilize naturally. Furthermore, due to Milnor Triality in AP Theory, (\cite{W1},p. 227), it is natural to refer to Artin presentations of the trivial group as 'bosons' and to Artin presentations of the binary icosahedral group of order $120$, as 'fermions'. Then, the SUSY which equates them, is false, but abstractly just barely so. Via Covering theory, then every fermion determines a boson. 

 According to Pauli, \cite{P}, p.ix, {\it "The Exclusion Principle is a necessary manifestation of the symmetries of the material world"}. The sheer existence of the discrete, purely group-theoretic, cone-like SUSY in AP Theory augments this view.

As intuited correctly in \cite{G}, SUSY is a consequence of low-dimensionality and group theory.

This is rigorously  butressed by the following analogues to the very important properties of classical SUSY mentioned in \cite{A} and \cite{Wi1}.

i) Relation with Dark Matter. See next section.

ii) AP Theory and Confinement (see section 5). The cone-like SUSY of AP Theory, just as that of Seiberg and Witten, see \cite{A}, is related to QCD confinement, in a topological analogous way.

iii) AP Theory and Donaldson/Seiberg-Witten Theory  (see p. 9 of \cite{W2})

iv) A promising topological relation with the  Standard Model \'a la \cite{BMS}, of Loop Quantum Gravity, see section 6, remark 2 below.

v) AP Theory is also a 'as holographic as possible', back-ground independent $(3+1)$-String Theory, (see \cite{Wi1}, p.3): Artin presentations as 'strings', which are inmune to the just nomographical criticisms given to mere continuous arcs in a back-ground space, as 'strings'. As presentations, the $r$ have a very intuitive string-like property caused by Jordan's Theorem in the plane, see \cite{C1}, p.364.

In AP Theory the discrete, group-theoretic 'string' $r$ is immediately related to $4D$ gravity, i.e., smooth $4D$ manifolds, via the the basic holographic construction of $W^4(r)$ from $r$, justifying the name {\it cosmic} strings. A priori, no higher dimensions are needed to give consistency to this 'string theory'.

As far as AP Theory is concerned, 't Hooft intuited correctly about the existence of a relation between $4D$ Quantum Gauge Theory and String Theory; see \cite{W2}, pp. 5, 11.

In AP Theory, the Hierarchy Problem ('why gravity is so weak', \cite{W}), can be solved, a priori, by observing, that the $(3+1)$-dimensional smooth space-time universes $W^4(r)$ are actually constructed holographically and non-locally, from the $2D$ membranes, the vacuum fluctuations $h(r):\Omega_n\to \Omega_n$. 

 The fact that AP Theory is cone-like, i.e. graded with a starting point, allows one to substitute mathematically rigorous  order and induction arguments, instead of non-rigorous {\it 'anthropic'} arguments, compare to \cite{Wi1}, p.6.

With AP Theory, in the cone-like String Landscape of the $W^4(r)$ universes, the non-perturbative, dynamic vacuum energy unleashed by the $\infty$-generated at each stage, cone-like graded group of Torelli topology- changing transitions/interactions, is fixed by the fundamental theory, without the need of an 'adjustment mechanism'. Compare to \cite{Po}.

\section{ Dynamic Dark Energy and the Cosmological Constant Problem}

Recall that in AP Theory, from the Artin presentations $r$ and the vacuum fluctuations/excitations, the $h(r):\Omega_n\to \Omega_n$, the smooth $4D$ universes $W^4(r)$, the last vestiges of Einsteinian gravity,
are constructed, in such a manner that the infinitely generated Torelli transitions/interactions act on all three {\it simultaneously, in unison} and, on the smooth $W^4(r)$, cause a radical and subtle homology-preserving $(3+1)$-topology-change. Thus we can say the Torelli transitions/interactions act on the holographic constituitivity of the $(3+1)$ space time universes $W^4(r)$. Compare to Schwinger's Newton quote above.

It is natural to consider this phenomenon as a rigorous qualitative analogue of Dynamic Dark Energy, the so-called energy of the vacuum, represented by a positive 'cosmological constant' and  non-rigorously atributed to the Uncertainty Principle, (see, e.g., \cite{Ca}, section $4$).

In other words, these {\it 'are good reasons to consider dynamic dark energy as an alternative to an honest cc',} see Carroll, \cite{Ca2}, p.11.

The Artin presentations $r$ form the discrete AP-vacuum and their, {\it infinitely generated at each stage},  group-theoretic Torelli symmetries give the dynamic energy of the vacuum, i.e., dynamic dark energy all the way  up to the infinitely generated at each stage, topology-changing transitions/interactions of the smooth $(3+1)$-universes, the $W^4(r)$.

Thus, in AP Theory, not only is {\it 'the quantum vacuum a dynamic medium'}, \cite{Wil1}, p.864, but, in AP Theory, the purely discrete Artin presentations, i.e. the 'empty space' in AP Theory, generates, via the $\infty$-generated at each stage, graded group of Torelli transitions, the dynamics of Dark Energy. 

This {\it qualitative} cc is a consequence of the cone-likeness, with topological Planckian starting point, of the extraordinary symmetry of AP Theory, the rigorous AP-analogue of classical 'broken' SUSY, given by the graded (infinitely generated at each stage) Torelli transitions/interactions. 

What other rigorous mathematical analogue can there be for the $(3+1)$-topology-change energy of the Torelli transitions/excitations acting on the purely discrete vacuum of the Artin presentations $r$, simultaneously, in unison, to also give smooth $(3+1)$-topology-change?

This qualitative meta-mathematical phenomenon is the 'field' analogue of the quantitative 'the cc is positive'. (\cite{Ca}, section 4) 

{\it The Artin Equation is the cone-like graded equation for vacuum energy.}

 Its group-theoretic, non-infinitesimal, non-analytic nature should be {\it necessary} for all the properties above to be possible and rigorously deducible. It is very unlikely that any non-trivial {\it metric} analytic differential-geometric PDE based theory could handle such radical smooth topology-changing transitions/interactions in unison with a topological Planckian starting point.

It is the sheer meta-mathematical existence of AP Theory, with its cone-like SUSY, its 'super gauge-ness' which is a qualitative, non-infinitesimal, {\it macroscopic} analogue for the positivity of the cosmological constant, as a measure of the energy density of the discrete group-theoretic vacuum, \cite{Ca}, p.7.

AP Theory is {\it 'the dynamical field'} which replaces a constant {\it 'fine-tuned'} numerical cc, p.33, \cite{Ca}.

Since unbroken classical Supersymmetry implies that the cc is zero, the sheer existence of AP Theory seems to indicate that also a classical rigorous  SUSY would have to be 'broken', but in a nice, canonical, graded with a starting point, oversee-able way by AP Theory: 'broken' is just its graded 'cone-like-ness'. Compare to Witten, \cite{Wi1}. Thus there is no need for non-rigorous 'anthropic' arguments here.

There is {\it 'no special symmetry which could enforce a vanishing vacuum energy, while remaining consistent with the known laws of physics'}, p.9, \cite{Ca}.

{\bf The only 'symmetry', which comes closest in doing this rigorously, is the sheer meta-mathematical existence of the very symmetric group-theoretic AP Theory itself as a whole.}

This again supports AP Theory as the incarnation, embodiment,  of cone-like SUSY. In AP Theory the 'positivity of the cc' is due to the graded cone-like-ness of SUSY and the power of its $\infty$-generated at each stage group of topology-changing transitions/interactions.

It is natural to suppose that as the sharpest $(3+1)$-QFT of section 2 above, with the most general and radical interactions, should also be the QFT of Dynamic Dark Energy.

 There is no need to {\it  postulate}  a {\it 'Quantum Higgs Field'}, \cite{W}, p.96, nor a {\it 'Quintessence Field'}, to start solving the conundrum of the cc problem. The whole AP Theory is that meta-mathematical 'field'.

It is interesting to ask whether the Poincar\'e homology 3-sphere, (as a 'composite graviton', \cite{Su}), whose fundamental group is the binary icosahedral group, $I(120)$, of order 120, is the cause, (when using the usual analytic field-theoretic methods), that the classic quantitative cc is 120 orders of magnitude larger than the observed one. Compare to p.31, \cite{Ca}.

 Notice the general semantic similarities of this rigorous, intrinsic, dynamic qualitative section $4$ with the more heuristic, quantitative, analytic approaches to the cc, e.g., \cite{BP}: {\it "... the repeated nucleation of membranes dynamically generates regions with a cosmological constant..."}.

\section{ The YM Millenium Problem with Mass Gap}

This section augments section 3 of \cite{W2}, by relating the Mass Gap problem with the cc problem by stressing the Dynamics of the Torelli transitions/interactions, as in section 4 above.

Due to its basic, conceptually simple, purely discrete group-theoretic nature, it is clear that AP Theory is a {\it model-free}  smooth topological $4D$ Quantum Yang-Mills Theory, in the sense of \cite{JW}. Due to all the properties above, AP Theory indeed gives a {\it 'mathematical understanding of the quantum behaviour of four-dimensional gauge theory, and a precise definition of quantum gauge theory in four dimensions}. Compare to \cite{JW}, p.3.

AP Theory contains, in particular, a {\it rigorous, non-perturbative, parameter-free, model-free} Quantum YM Gauge Theory, which has cone-like SUSY, making it even more of a incarnation of symmetry than QCD, as in \cite{W}, p.196.

AP Theory is indeed a mathematically complete example of a smooth topological, model-independent quantum gauge theory in four-dimensional spacetime, a precise definition of a sporadic, intrinsic, model-free, quantum gauge theory characteristic of four dimensions.

Meta-mathematically, this already solves topologically, the first part of the YM Millenium problem, without any 'new fundamental ideas'. Compare to \cite{JW2}. The ideas of Poincar\'e, naturally augmented by AP Theory, suffice. See section 6 ahead.

 The purely group-theoretic, non-analytic Artin Equation avoids having to solve intractable PDE and other nomographical problems of QCD in its present state and in other theories. Its discreteness immediately gives non-perturbativity as well as back-ground independence and the powerful topology-changing cone-like, graded group of Torelli transitions/interactions, where $\infty$ generation, at each stage, replaces the rigid, static $\infty$ of the dimension of classical Hilbert space.

The serious problem of {\it 'the masslessness of classical YM-waves'}, \cite{JW}, p.2, does not arise AP Theory: the membrane-ness of the basic topological Planck particle-membrane duality gives 'mass', to these 'waves', wavicles', the $h(r):\Omega_n\to \Omega_n$, the membranic vacuum fluctuations.

Thus it is natural to expect AP Theory to be the theory in which the YM Millenium Mass Gap  problem should be solved, at least topologically, qualitatively, {\it dynamically} in a model-free manner, in analogy to section 4 above.

We now point out, building on section 4 of \cite{W2}, how the second part of the YM Millenium Problem, the so-called 'mass gap problem', also has rigorous qualitative, {\it dynamic} topological meaning, intimately related to the cc.

In AP Theory the YM Millenium Mass Gap problem is the $2D$ meta-mathematically analogue, parallel to the $4D$ cc problem of section 4.

The mass gap is given qualitatively by the infinitely generated, at each stage, Torelli transitions/interactions, (whose dynamic $\infty$ substitutes the rigid static $\infty$ of the dimension of Hilbert space), acting on the vacuum fluctuations $h(r):\Omega_n\to \Omega_n$ in unison with their dynamic dark energy action on the spacetimes $W^4(r)$.

 We view AP Theory more 'locally', (just with the $h(r):\Omega_n\to \Omega_n$), as a {\it rigorous} discrete group-theoretic analogue of QCD, which, as already mentioned,  (see p.196 of \cite{W}), {\it "..is in a profound and literal sense constructed as the embodiment of symmetry"}.

The cause of a Mass Gap in AP Theory is just the dynamic $2D$ energy displayed by the the Torelli transitions, acting  on the vacuum fluctuations $h(r):\Omega_n\to \Omega_n$, instead of the space-time universes $W^4(r)$, where they generate dynamic dark energy.

Qualitatively, dynamically, topologically speaking, a positive mass gap is related to a positive cosmological constant, via the Torelli action acting simultaneously, in unison on the 'local' vacuum fluctuations, the $h(r):\Omega_n\to \Omega_n$ and the 'global' spacetime universes, the $W^4(r)$. Compare to \cite{J}.

 The Torelli just acting on the $h(r):\Omega_n\to \Omega_n$, i.e.,  acting membranically, give "massive excitations" in unison,  while generating Dark Energy (topology change) on  the space-time universes $W^4(r)$.

This is the {\it qualitative} analogue to the {\it quantitative} mass gap of property (1) of p.3 of \cite{JW}.

{\bf In AP Theory, considered as a Quantum Yang-Mills Theory, the Mass Gap is intimately related to Dynamic Dark Energy, in a very concise non-perturbative  manner and the mathematical reason for their dynamic 'positivity', via the powerful Torelli action is essentially the same.}

The cc problem and the YM Millenium Mass Gap problem have exactly the same {\it dynamic} cause and this is only made possible by substituting the {\it static} rigidity of the $\infty$ dimension of Hilbert space, by the graded, cone-like {\it dynamic} $\infty$ generation, at each stage, of the group of Torelli transitions/interactions.

Thus the mass gap in AP Theory is 'dynamically generated' as in more classic analytic perturbative approaches, see \cite{Go}, \cite{F}.

The mass/energy of these 'massive excitations' is stronger than merely letting the $h(r):\Omega_n\to \Omega_n$
iterate themselves.

All this so that genuine topological analogues to other important properties, besides the mass gap, (i.e.,property (1) of p.3 of \cite{JW}) of the YM Millenium problem, properties (2) and (3) also mentioned on p.3 of \cite{JW}, are displayed as follows. 

Property (3) is trivial in AP Theory, the Torelli obviously leaving the cone-like graded vacuum invariant.

Property (2), {\it confinement}, seems to be the clear also, as explained in \cite{W2}, section 3: due to the fact that the $h(r):\Omega_n\to\ \Omega_n$ are determined only up to isotopies, their fixed points, i.e. our 'quarks', 'partons', although they have to exist mathematically, can not be determined individually, i.e., they are invisible so to speak. See \cite{JW2}.

{\it In AP Theory 'confinement'  follows immediately from 'indeterminism'}

 The vacuum fluctuations/ excitations $h(r):\Omega_n\to \Omega_n$, when iterated, are topological versions of the  Casimir Effect.  Compare to \cite{F}, \cite{J}.

Thus in AP Theory, the Casimir effect, the mass gap and dynamic dark energy are all related via the Torelli transitions/interactions, operating on the membranic vacuum fluctuations $h(r):\Omega_n\to \Omega_n$.

This scenario becomes even more dynamically chaotic, since these Torelli actions combine naturally with the mere iteration of the homeomorphisms $h(r):\Omega_n\to \Omega_n$, the AP analogue of the Casimir force.

Thus Nielsen-Thurston theory as well as Artin-Mazur zeta functions come into play.

In this section we have shown that inside AP Theory lies a rigorous conceptually simple topological theory with topological properties, which are genuine qualitative analogues to the most important desiderata of the actual, analytic, quantitative YM Millenium Problem.

The Mass Gap part of the YM Millenium problem is nothing but the qualitative 'positive cc' part of the Dynamic Dark Energy problem as in section 4, but at the $2D$ level of the membranic vacuum fluctuations $h(r):\Omega_n\to \Omega_n$, instead of at the $4D$ level of smooth topology-change of the space-time universes $W^4(r)$.

Given the fact that, as explained in the introduction, the rigidity of the $\infty$ dimension of classical Hilbert space, (which is used in the formulation of the classical YM Millenium Mass Gap problem, (see, e.g., \cite{V}) has been dynamically broken by the graded, $\infty$-generated at each stage, Torelli transitions/interactions, this {\it dynamic} solution to the mass gap problem is the best one can hope for in AP Theory.

 Due to its intimate meta-mathematical similarity wih AP Theory as a Dynamic Dark Energy Theory, due to the holographic construction of the $(3+1)$-universes $W^4(r)$ from the mass excitations of the vacuum transitions/excitations, the $h(r):\Omega_n\to \Omega_n$, and the fact that the Torelli actions act in unison on both, it seems unlikely that the original YM Millenium Problem of \cite{JW} can be completely solved analytically, in a rigorous {\it non-perturbative} manner.

It is very likely that the discreteness of the purely group-theoretic Artin Equation as well as the substitution in AP Theory, of the rigid $\infty$ of the dimension of classical Hilbert space, by the cone-like graded dynamic $\infty$ of the generation, at each stage, of the topology-changing Torelli transitions/interactions, are  necessary for all these properties to co-exist rigorously, even just in the topological case, leave alone the more rigid {\it metric} cases. Compare to \cite{W}, p.97. Compare also to remark 1 below.

{\it Does AP Theory destroy, with its infinitely generated, smooth, topology-changing transitions, any hope of finding smooth analytic, non-perturbative solutions to any physically relevant equations appearing in the actual YM Millenium problem?}  

In other words, it seems very likely that the starting point of a rigorous, classic YM Millenium Theory  of \cite{JW}, has to reside in AP Theory (and be related there to Dynamic Dark Energy Theory), but once inside it, it would acquiere the seeds of its own destruction, {\it by the $\infty$-generated topology-changing Torelli dynamics}, as a mathematically rigorous, analytic, smooth non-perturbative $(3+1)$-QFT.

{\it We conjecture that any non-perturbative, smooth analytic approach to the $(3+1)$-YM Millenium problem is incompatible with the universal, $\infty$-generated at each stage, graded smooth topology-change theory of the Torelli transitions/interactions}

{\it Does the AP topological Planck analogue, i.e., the AP Theory, $r$/ $h(r):\Omega_n\to \Omega_n$, 'cosmic string/membrane' duality, 'do' for the cc problem, (as in \cite{Ca}), and the YM Millenium Mass Gap problem, (as stated in \cite{JW}), what Planck's Law 'did' for the Raleigh-Jeans UV problem?}

\section{Historical Remarks, Questions, Conjectures}.

1. Due to the basic conceptual simplicity of AP Theory, all of whose concepts were known very well to the founder of modern topology, Poincar\'e, it is perhaps not suprising that his semi-abstract 'resonateurs', (oscillators), which he never defined explicitly, (compare to \cite{P1}, p.654, \cite{Pl}, pp.388, 390, \cite{I}, \cite{Pr}, p.340) seem to be topologically realized by the vacuum fluctuations/excitations, the $h(r):\Omega_n\to \Omega_n$. All this without having to introduce the rigidity of infinite-dimensional Hilbert space of later Quantum Mechanics. In this context, the non-analytic, purely discrete group-theoretic nature of the fundamental Artin Equation corresponds to the proof by Poincar\'e (\cite{P1}, pp. 678-681, \cite{Pr}) and Ehrenfest of the actual theoretical necessity of Planck's original {\it postulate}, relating continuity and discreteness.

Just like Poincar\'e does not need the main concepts of thermodynamics, see Planck, \cite{Pl}, \cite{Pr}, p.340, \cite{I}, p.880, neither does AP Theory.
 
If one also notes that Poincar\'e, the founder of the topological theory of manifolds, who moreover introduced discrete group theory therein, was fascinated with early Quantum Theory, one almost feels tempted to say "Poincar\'e missed an opportunity", as far as discovering the existence of AP Theory is concerned. 

AP Theory could be called "Poincar\'e's Smooth, Holographic $(3+1)$-TOE", since he invented and/or introduced literally all of its basic mathematical concepts. 

2. Rigorous, non-perturbative purely group-theoretic calculations can be done in AP Theory, \cite{W1}, \cite{W3},  \cite{CW}, \cite{C}. Thus AP Theory and its pure framed/colored braids should be related to so-called {\it 'preons, constituted of braids of space-time'} which {\it 'survive quantum fluctuations'}, etc. of Loop Quantum Gravity; compare to \cite{BMS}. 

Their Standard Model versions would be immediately related to $(3+1)$-gravity via the smooth, $(3+1)$-spacetimes $W^4(r)$, thus giving a $3+1)$-TOE.
 In AP Theory $(3+1)$-gravity seems feeble, weak, (\cite{W}, p.151) because of its  holographic nature, since $4D$ smooth structures of the $W^4(r)$ are in function of $2D$ smooth structures, of the membranic vacuum fluctuations/excitations, the $h(r):\Omega_n\to \Omega_n$. This seems relevant to the so-called Hierarchy Problem, as well as 'mass without mass' arguments, (\cite{W}, p.129), in QCD.

3. AP Theory has many topological gauge-theoretic similarities with the Quantum Hall Effect (QHE), the vacuum fluctuations $h(r):\Omega_n\to \Omega_n$ serving as the 'plateaus', (see \cite{PG}, pp. 17, 18, 103, 355; \cite{BES}, p.5374 and \cite{MM}), thus extending Laughlin's geometric gauge theory of the QHE  in a universal, intrinsic, non-perturbative, model-free, gauge-theoretical topological manner.

In this context Calcut's thesis, \cite{C}, can be interpreted as a universal, global, intrinsic {\it gravitational} QHE.

4. Although one can compute exactly in a finite, non-perturbative manner in AP Theory,using, e.g., the computer algebra system {\bf magma}, (\cite{W1}, \cite{W2}, \cite{W3}, \cite{C}, \cite{C1}, \cite{CW}), the only {\it quantitative} prediction AP Theory makes at this point is: if the Higgs boson appears experimentally, it will be somehow related to Poincar\'e's homology 3-sphere and representations of its fundamental group, the binary icosahedral group, $I(120)$, of order 120.


\begin{thebibliography}{WWWWWW}

\bibitem[A]{A} Ananthaswamy, A., \emph{In SUSY we trust:What the LHC is really looking for}, New Scientist, 11 November, 2009. \smallskip

\bibitem[ACK]{ACK} Arkadi-Hamed, N., Cachazo, F., Kaplan, J., \emph{What is the simplest Quantum Field Theory?}, arXiv:0808.1446v2 [hep-th] 19 Aug 2008. \smallskip

\bibitem[BES]{BES} Bellisard, J., v.Elst, A., Schulze-Baldes, H., \emph{The non-commutative geometry of the Quantum Hall Effect}, J.Math.Phys.\text{35}, 5374-5451. \smallskip

\bibitem[BJSV]{BJSV} Bershadsky, M., Johansen, A., Sadov, C., Vafa, C. , \emph{Topological Reduction of $4D$ SYM to $2D$ $\sigma$-Models}, Nucl. Phys. B448 (1995) 166-186. hep-th/9501096. \smallskip

\bibitem[BMS]{BMS} Bilson-Thompson, S., Markopolou, F., Smolin, L., \emph{Quantum Gravity and the Standard Model}, hep-th/0603022. \smallskip

\bibitem[BP]{BP} Bousso, R., Polchinski, P., \emph{Quantization of four form fluxes and Dynamical Neutralization of the cosmological constant}, hep-th/0004134. \smallskip


\bibitem [C]{C} Calcut, J. S., \emph{Torelli Actions and Smooth Structures on 4-manifolds}, J.Knot Theory Ramifications \text{17}, (2008), 171-190. \smallskip

\bibitem[C1]{C1} Calcut, J. S., \emph{Artin Presentations from an Algebraic Viewpoint}, Journal of Algebra and Its Applications, Vol. 6, No.2 (2007) 355-367. \smallskip 

\bibitem[CFW]{CFW} Calegari, D. , Freedman, M., Walker, K. , \emph{Positivity of the Universal Pairings in 3 Dimensions.} J.Amer. Math. Soc. \text{23} (2010), 107-188. \smallskip

\bibitem[Ca]{Ca} Carroll, S., \emph{The Cosmological Constant}, Living Rev. Relativity, \text{3}, 2001.  \smallskip

\bibitem[Ca2]{Ca2} Carroll, S., \emph{Why is the Universe Accelerating?}, In Measuring and Modeling the Universe, 2004, W. L. Freedman,ed., Cambridge University Press. \smallskip


\bibitem[CW]{CW} Calcut, J. S. and Winkelnkemper, H. E., \emph{Artin Presentations of Complex Surfaces}, Bol. Soc. Mat. Mexicana \text{10}, Special issue in honor of F. Gonz\'alez-Acu\~na, 2004, 63-87. \smallskip

\bibitem[De]{De} DeWitt, B., \emph{Quantum Gravity, Yesterday and Today}, arXiv:0805.2935v1 [physics.hist-ph] 19 May 2008. \smallskip

\bibitem[D]{D} Dyson,F., \emph{Birds and Frogs} Not. AMS. \text {56}, no.2, 2009, 212-223. \smallskip

\bibitem[Fe]{Fe} Feshbach, H., \emph{Group Theory:Who needs it?}, Physics Today, August 1986. \smallskip

\bibitem[F]{F} Faddeev, L., \emph{Mass in quantum Yang-Mills theory (comment on a Clay Millenium Problem)}, Perspectives in analysis,63-72, Math.Phys.Stud.,27, Springer, Berlin, 2005. \smallskip

\bibitem[G]{G} Gates. Jr., S. J. , \emph{On the Possible Group Theoretic and Low-dimensional Origins of Spacetime Supersymmetry}, Talk at MSRI, Berkeley, April 2002. \smallskip

\bibitem[Go]{Go} Gogokhia, V., \emph{The Mass Gap Problem and Solution to the Gluon Confinement Problem in QCD}, arXiv:hep-ph/0702066v2 22 Feb 2007. \smallskip

\bibitem[H]{H} 't Hooft, G., \emph{The mathematical basis for deterministic quantum mechanics}, arXiv:quant-ph/0604008v2 26 Jun 2006. \smallskip

\bibitem[I]{I} Irons, F. E., \emph{Poincar\'e's proof of quantum discontinuity interpreted as applying to atoms}, Am. J. Phys. \text{69}, 2001, 879-884. \smallskip

\bibitem[J]{J} Jaffe, R., \emph{The Casimir Effect and the Quantum Vacuum}, arXiv:hep-th/0503158. \smallskip

\bibitem[JW]{JW} Jaffe, A., Witten, E., \emph{Quantum Yang-Mills Theory}, preprint, 2004. \smallskip

\bibitem[JW2]{JW2} Jaffe, A., Witten, E., \emph{Yang-Mills Theory}, Mathematical Description of the Clay YM Millenium Problem, 1 page. \smallskip

\bibitem[LM]{LM} Leininger, C., Margalit, D., \emph{Two generator subgroups of the pure braid group}, to appear in {\it Geometria Dedicata.} \smallskip

\bibitem[MO]{MO} Montonen, C., Olive, D., \emph{Magnetic Monopoles are Gauge Particles?}, Ref.TH.2391-CERN 1977. \smallskip

\bibitem[MM]{MM} Marcolli, M., Mathai, V., \emph{Orbifold Fractional Quantum Hall Effect}, preprint, 2003. \smallskip

\bibitem[NN]{NN} Nielsen, H., Ninomiya, M., \emph{Card Game restriction in LHC can only be successful}. arXiv:0910.0359v1 [physics.gen-ph] 2 Oct 2009. \smallskip

\bibitem[P1]{P1} Poincar\'e, H., \emph{Sur la theorie des Quanta}, Oeuvres, t.9, 1912, 626-623. \smallskip

\bibitem[P2]{P2} Poincar\'e, H., \emph{Sur la Dynamique de l' Electron}, Oeuvres, t.9, 1905, 494-550. \smallskip

\bibitem[Po]{Po} Polchinski, J., \emph{The Cosmological Constant and the String Landscape}. arXiv:hep-th/0603249v2 21 Apr 2006. \smallskip

\bibitem[P]{P} Pauli, W., \emph{Collected Scientific Papers}, v.1, R. Kronig, V. F. Weisskopf, eds. \smallskip

\bibitem[Pl]{Pl} Planck, M., \emph{Poincar\'e und die Quantentheorie}, Acta Math. \text{38}, 1921, 387-397.\smallskip

\bibitem[PG]{PG} Prange, R., Givrin, S., eds., \emph{The Quantum Hall Effect}, Graduate Texts in Contemporary Physics, Springer, 1987. \smallskip

\bibitem[Pr]{Pr} Prentiss, J., \emph{Poincar\'e's proof of the quantum discontinuity of Nature}, Am. J. Phys. \text{63}, 1995, 339-350. \smallskip

\bibitem[Su]{Su} Sundrum, R., \emph{Towards an effective particle-string resolution of the cosmological constant problem}, J.High Energy Phys. \text{1999}, hep-ph/9708329. \smallskip

\bibitem[S]{S} Schwinger, J., \emph{A Magnetic Model of Matter}, Science \text{165}, 22 August 1969, 757-761. \smallskip

\bibitem[V]{V} Vogt, E., \emph{Existenz von Quanten-Yang-Mills Theorien mit Massenluecke} Elem.Math. \text{57}, 2002, 121-126. \smallskip

\bibitem[We]{We} Weyl,H., \emph{Gesammelte Abhandlungen}, IV. \smallskip

\bibitem[W]{W} Wilczek, F., \emph{The Lightness of Being}, Basic Books, 2008. \smallskip

\bibitem[Wil1]{Wil1} Wilczek, F., \emph{Asymptotic Freedom: From Paradox to Paradigm}, Nobel Lecture, Rev.Mod. Phys., Vol 77, No.3, July 2005. \smallskip

\bibitem[Wil2]{Wil2} Wilczek, F., \emph{The Origin of Mass}, MIT Physics Annual 2003. \smallskip

\bibitem[W1]{W1} Winkelnkemper, H. E., \emph{Artin Presentations,I: Gauge Theory, $(3+1)$ TQFTs and the Braid Groups}, J.Knot Theory Ramifications \text{11}, (2002), 223-275. \smallskip

\bibitem[W2]{W2} Winkelnkemper, H. E., \emph{AP Theory II: Intrinsic 4D Quantum YM Theory with Mass Gap}, arXiv:0711.2054v2 [math.GT] 22 Feb 2008. \smallskip

\bibitem[W3]{W3} Winkelnkemper, H. E., \emph{What is...an Artin Presentation?}, preprint, 2003. Available at http//www.math.umd.edu/\~ hew/. \smallskip

\bibitem[Wi1]{Wi1} Witten, E., \emph{Supersymmetry and other scenarios}, preprint. \smallskip

\bibitem[Wi2]{Wi2} Witten, E., \emph{Is Supersymmetry Really Broken?}, hep-th/9409111. \smallskip

\bibitem[Wi3]{Wi3} Witten, E., \emph{Next Steps in Particle Physics}, preprint, 2009. \smallskip

\bibitem[Wi4]{Wi4} Witten, E., \emph{String Theory}, preprint, October 2001. \smallskip

\bibitem[Wi5]{Wi5} Witten, E., \emph{Gauge Theory and the Geometric Langlands Program}, preprint, 2005. \smallskip





 



\end{thebibliography}
\end{document}